\documentclass[twoside,12pt]{article}

\usepackage[all]{xy}
\usepackage{amssymb,amsmath,latexsym}
\usepackage{theorem}
\usepackage{amsfonts}
\usepackage{ulem}
\usepackage{amscd}
\usepackage{amsxtra}
\usepackage{mathrsfs}

\newtheorem{thm}{Theorem}[section]
\newtheorem{cor}[thm]{Corollary}
\newtheorem{lem}[thm]{Lemma}
\newtheorem{prop}[thm]{Proposition}

\newtheorem{rem}[thm]{Remark}

\numberwithin{equation}{section}

\newcommand{\cqfd}
{\hspace{1cm}
\rule{2mm}{2mm}%
\medbreak%
\par%
}

\def\pr{{\parindent0pt {\bf Proof.\ }}}
\def\cqfd
{\hspace{1cm}
\rule{2mm}{2mm}%
\medbreak%
\par%
}

\author{}

\begin{document}
\title{Dhara-Rehman-Raza's identities on left ideals of prime rings}

\date{}
\maketitle \vspace*{-1.5cm}

\thispagestyle{empty}


\begin{center}
\author{Driss Bennis$^{1}$, Brahim Fahid$^{2}$  and  Abdellah Mamouni$^{3}$}

\begin{center}
\small{1:  Centre de Recherche de Math\'ematiques et Applications de Rabat (CeReMAR), Faculty of Sciences, Mohammed V University in Rabat,  Morocco.}
\end{center}\vspace{-0,4cm}
\hspace{2cm} \small{ d.bennis@fsr.ac.ma; driss$\_$bennis@hotmail.com}\medskip

\small{2:   Department of Mathematics, Superior School of Technology, Ibn Tofail University, Kenitra, Morocco.}\\
 $\mbox{}$ \small{fahid.brahim@yahoo.fr}\medskip

\small{3:   Department of Mathematics, Faculty of Sciences and Techniques, Moulay Isma\"il  University, Errachidia, Morocco.}\end{center}\vspace{-0,4cm}
$\mbox{}$ \hspace{4cm}  \small{mamouni\_1975@live.fr}\smallskip

\bigskip\bigskip
\noindent{\large\bf Abstract.}
It is known that every nonzero Jordan ideal of $2$-torsion free semiprime rings contains a nonzero ideal. In this paper we show that also any square closed Lie ideal  of  a $2$-torsion free prime ring   contains a nonzero ideal. This can be interpreted by saying that studying identities over one sided ideals is the ``optimal" case to study identities. With this fact in mind, we generalize some results of Dhara, Rehman  and   Raza in [Lie ideals and action of generalized derivations in rings, Miskolc Mathematical Notes, \textbf{16} (2015), 769 -- 779] to the context of nonzero left ideals. \bigskip


\section{Introduction}
Let $R$ be an associative ring with center $Z(R)$. Recall that $R$ is prime if for $a, b \in  R$, $aRb = \{0\}$ implies either $a = 0$ or $b = 0$.   $R$ is said to be semiprime if for $a \in R$, $aRa =\{0\}$ implies $a = 0$. Let $R$ be a prime ring. For any pair of elements $x, y \in  R$, we shall write $[x,y]$ (resp., $x \circ y $) for the commutator $xy-yx $ (resp., for the Jordan product $ xy+yx $). An additive subgroup $J$ of $R$ is said to be a Jordan ideal of $R$ if $ u \circ r \in  J$ for all $ u \in J$ and $r \in R$. Note that every ideal of $R$ is a Jordan ideal of $R$ but the converse is not true in general. An additive subgroup $U$ of $R$ is said to be a Lie ideal of $R$ if $[u,r] \in U$ for all $u \in U $ and $ r \in R$. It is clear that if characteristic of $R$ is $2$, then Jordan ideals and Lie ideals of $R$ coincide. A Lie ideal $U$ of $R$ is said to be square closed if $u^2 \in U $ for all $u \in U$.\medskip

An additive mapping $ d : R \longrightarrow R $ is called a derivation if $d(xy) = d(x)y + xd(y)$ holds for all $ x,y \in R$. In \cite{B} Bre\v{s}ar introduced the  generalized derivation: an additive mapping $ F : R \longrightarrow R $ is called a generalized derivation if there exists a derivation $ d : R \longrightarrow R $ (an associated derivation of $F$) such that $ F(xy) = F(x)y + xd(y)$   for all $ x, y \in R$.    The notion  a generalized derivation covers both the notions of a derivation and of a left multiplier (i.e., an additive mapping $ f : R \longrightarrow R $ satisfying $f(xy) = f(x)y$ for all $x, y \in R$). A ring $R$ is said to be $n$-torsion free, where $ n \neq 0$ is a positive integer, if whenever $na = 0$, with $a \in R$, then $ a = 0$.\medskip

In the present paper we investigate commutativity of a prime ring satisfying certain differential identities on a nonzero left ideal. Let us first recall that  the study of commutativity of rings using differential identities goes back to the well-known Posner's work \cite{P57}  in which he proved that the existence of a nonzero centralizing derivation $f$ on a prime ring $R$ (i.e.,   $[f(x),x]\in Z(R)$ for all $x\in R$) forces the ring $R$ to be commutative. This result is now known by Posner's second theorem. Since then, several authors have been interested in   extending or  generalizing this result to  different contexts. For instance,  Awtar, in \cite{A}, extended Posner's second theorem to the case of centralizing derivation on either nonzero Jordan ideals or nonzero square closed Lie ideals (see also Mayne's papers \cite{M1,M2,M3}).  In \cite{B93}, Bre$\mathrm{\check{s}}$ar generalised Posner's second theorem to identities related to two derivations (see  \cite[Theorem 4.1]{B93}).  In \cite{H98}, Hvala introduced the study of identities related to generalized derivations. In \cite{AFS}, Ali, De Filippis and Shujat studied identities with generalized derivations  on one sided ideals of a semiprime ring.\medskip

Also motivated by the success that known Posner's second theorem, several authors have introduced new kind of  differential identities. In this context,  Ashraf and Rehman proved in \cite{AR} that a prime ring $R$ with a nonzero ideal $I$ must be commutative if $R$ admits a nonzero derivation $d$ satisfying $d(xy)-xy \in
Z(R)$ for all $ x,y \in  I$ or $d(xy)+xy \in Z(R)$ for all $x,y \in  I$. In \cite{AAA}, Ashraf, A. Ali    and S. Ali studied these identities in the case where $d$ is a generalized derivation. In \cite{DRR},  Dhara, Rehman and Raza used  more general differential identities. Precisely, they    showed that for a nonzero square closed Lie ideal $U$ of a prime ring $R$, if $R$ admits  nonzero generalized derivations $F$, $G$ and $H$  satisfying $F(x)G(y)\pm H(xy)\in Z(R)$ or $F(x)F(y)\pm H(yx)\in Z(R)$ for all $x, y\in U$, then $U \subseteq  Z(R)$. See also      \cite{AliDharaKhan,   DKM, DKP, DPA,    TSD} for other works related to these identities.\medskip

Then, naturally one can ask whether we get the same conclusion as the one of Dhara, Rehman and Raza's results if we replace the nonzero square closed Lie ideal of the prime ring $R$ by  other particular subsets of $R$. Namely,  Jordan ideals, (both sided) ideals and one sided ideals. It is important to mention that   the study of  identities on  Jordan ideals and (both sided) ideals can be considered as particular case of the study  of these identities on square closed Lie ideals. Indeed,  it is clear that  every  ideal is a square closed Lie ideal and, by \cite[Theorem 1.1]{H}, every nonzero Jordan ideal of $2$-torsion free semiprime rings contains a nonzero ideal. Then, only the ``one sided ideals" case could be of interest. Our aim in this paper is to show that  $R$ will be commutative if we consider, in Dhara-Rehman-Raza's identities, only   nonzero left ideals instead of  square closed Lie ideals (see Theorems \ref{th-one} and \ref{th-tow}). This is a generalization of  Dhara,  Rehman and   Raza's results (see Corollary \ref{cor-princ3}) we prove  that  any square closed Lie ideal $L$ of  a $2$-torsion free prime ring $R$  contains a nonzero ideal of $R$ (see Proposition \ref{prop-2L}).

\section{Preliminaries}
In this section we recall and present some properties of  left ideals which will be used to prove our main results.

\begin{lem}\label{lem-p1}
    Let $R$ be a prime ring, $I$ is a nonzero left ideal of $R$ and  for every $a,b \in R$, $aIb = (0)$, then $a = 0$ or $Ib=(0)$.
\end{lem}
\pr
We have $aIb = (0)$ imply $aRIb = (0)$, so by primeness of $R$ we get $a = 0$ or $Ib=(0)$.
\cqfd
\begin{lem}\label{lem-p2}
    Let $R$ be a prime ring and $I$ a nonzero left ideal of $R$. If $F$ is a generalized derivation of $R$ with an associated derivation $d$ such that $F(I) = 0$, then $Id(I) = 0$.
\end{lem}
\pr
We have,  for every $u,v \in I$,  $F(uv) =F(u)v+ud(v)=ud(v)= 0$.  Hence $Id(I) = 0$.
\cqfd
\begin{lem}\label{lem-p3} Let $R$ be a prime ring and $I$ be a nonzero left ideal of $ R$.  Then, the following assertions are equivalents:
\begin{enumerate}
  \item $[I,I]=0$.
   \item $ I\subseteq Z(R)$.
   \item $R$ is commutative.
\end{enumerate}
\end{lem}
\pr
$(1)\Rightarrow (2)$ Let $I$ be a left ideal such that $[I,I]=0$. Then,  $0=[I,RI]=[I,R]I$, which implies that  $[I,R]=0$. \\
$(2)\Rightarrow (3)$  We have $[I,R]=0$. This gives $0=[RI,I]=[R,R]I$. Since a left annihilator of a left ideal is zero,  $[R,R]=0$;  that is, $R$ is commutative.\\
$(3)\Rightarrow (1)$ Obvious.
\cqfd
\begin{lem}\label{lem-p4}
  Let $R$ be a  noncommutative prime ring and $I$ be a nonzero left ideal of $ R$. Let $d$  be  derivation of $R$ and $z\in Z(R)$,  such that $[xy,r]d(z)=0$ for all $x, y \in I$ and $r\in R$. Then $d(z)=0$.
\end{lem}
\pr
First we show that $d(z)\in Z(R)$, for this we have $d(zr)=d(z)r+zd(r)=d(r)z+rd(z)$, so $d(z)\in Z(R)$. On other hand we have
\begin{equation}
  [xy,r]d(z)=0
   \;\;\;\mbox{for all}\;\;x, y \in I;\;\; r\in R.
\end{equation}
Replacing $x$ by $sxd(z)$ in the above equation, where $s \in R$, we get
\begin{equation}
[s,r]xd(z)yd(z)=0
   \;\;\;\mbox{for all}\;\; x, y \in I;\;\; r\in R.
\end{equation}
So $Id(z)=0$ and $d(z)=0$.
\cqfd

It is well-known that  every nonzero Jordan ideal of a $2$-torsion free semiprime ring contains a nonzero ideal (see, for instance, \cite[Theorem 1.1]{H}). Here we give a similar result for square closed Lie ideals.

\begin{prop}\label{prop-2L}
 Let $L$ be a nonzero square closed Lie ideal of a $2$-torsion free prime ring $R$. Then, $L$ contains a nonzero ideal of $R$.
\end{prop}
\pr
Let $L$ be a Lie ideal of $R$  such that $u^2 \in L$ for all $u \in L$. Therefore, for any $u,v \in L$, we get $uv+vu=(u+v)^2-u^2-v^2\in L$. On the other hand we have $uv-vu \in L$. Combining these two equalities we get $2uv \in L$ for all $u,v\in L$. Then,  $2L$ is both a Lie ideal and a subring of $R$. Also $2L\neq 0$ since $R$ is a $2$-torsion free ring. Then, by \cite[Lemma 1.3]{H} and \cite[Lemma 5]{DRR},  $2L$ contains a non-zero ideal of $R$ and so does $L$ since $2L\subseteq L$.  \cqfd


\section{Main results}

We start with the first main result, using Proposition \ref{prop-2L}, it can be seen as a generalization of  \cite[Theorems 1 and 2]{DRR} ( see Corollary \ref{cor-princ3}).

 \begin{thm}\label{th-one}  Let $R$ be a  prime ring and $I$ be a nonzero left ideal of $ R$. Let $F,$ $G$ and $H$ be generalized derivations associated to  derivations $f,$ $g$ and $h$ of $R$, respectively, such that $F(x)G(y)-H(xy)\in Z(R)$ for all $x, y\in I$.  Then $R$ is commutative or $If(I)= 0$ or $Ig(I)= 0$.
 \end{thm}
\pr
We are given that
\begin{equation}\label{eq1}
  F(x)G(y)-H(xy)\in Z(R)
   \;\;\;\mbox{for all}\;\;x, y \in I.
\end{equation}
 Replacing $y$ by $yz$ in (\ref{eq1}),  we get
\begin{equation}\label{eq2}
  F(x)G(y)z-H(xy)z+F(x)yg(z)-xyh(z)\in Z(R)
   \;\;\;\mbox{for all}\;\;x, y, z \in I.
\end{equation}
That is
\begin{equation}\label{eq3}
  [F(x)yg(z)-xyh(z),z]=0
   \;\;\;\mbox{for all}\;\;x, y, z \in I.
\end{equation}
Replacing $x$ by $xu$, where $u\in I$, we get
\begin{equation}\label{eq4}
  [F(x)(uy)g(z)-x(uy)h(z),z]+[xf(u)yg(z),z]=0
   \;\;\;\mbox{for all}\;\;u, x, y, z \in I.
\end{equation}
Since $uy \in I$, then (\ref{eq3}) together with (\ref{eq4}) force that
\begin{equation}\label{eq5}
  [xf(u)yg(z),z]=0
   \;\;\;\mbox{for all}\;\;u, x, y, z \in I.
\end{equation}
Writing $rx$ instead of $x$ in (\ref{eq5}), where $r\in R$, we obtain
\begin{equation}\label{eq6}
  [r,z]xf(u)yg(z)=0
   \;\;\;\mbox{for all}\;\;u, x, y, z \in I; \;\; r,s\in R.
\end{equation}
That is
\begin{equation}\label{eq7}
  [r,z]If(u)Ig(z)=0
   \;\;\;\mbox{for all}\;\;u, z \in I; \;\; r,s\in R.
\end{equation}
The primeness of $R$ together with (\ref{eq7}) imply that $[r,z]=0$ or $If(u)=(0)$ or $Ig(z)=(0)$. Therefore, $R$ is commutative or $If(I)= 0$ or $Ig(I)= 0$.
\cqfd

As a consequence of Theorem \ref{th-one} we get the following result.

\begin{cor}\label{cor-one}
 Let $R$ be a  prime ring and $I$ be   a nonzero left ideal of $ R$. Let  $F,$ $G$ and $H$ be generalized derivations associated to derivations $f,$ $g$ and $h$ of $R$,  respectively, such that $F(x)G(y)+H(xy)\in Z(R)$ for all $x, y\in I$. Then $R$ is commutative or $If(I)= 0$ or $Ig(I)= 0$.
\end{cor}
\pr
We notice that $-H$ is a generalized derivation  of $R$   associated to the  derivation $-h$. Hence replacing $H$ by $-H$ in Theorem \ref{th-one}, we get $ F(x)G(y)-(-H)(xy)\in Z(R)$ for all $x, y\in I$,  that is $ F(x)G(y)+H(xy)\in Z(R)$ all $x, y\in I$. Therfore  $R$ is commutative or $If(I)= 0$ or $Ig(I)= 0$.
\cqfd

Now we can see that \cite[Theorems 1 and 2]{DRR} can be seen as a consequence of Theorem \ref{th-one} as follows.

\begin{cor}\label{cor-princ3}
 Let   $L$ be a nonzero square closed  Lie ideal of  $2$-torsion free prime ring $R$. Let $F,G,H $ be   generalized derivations with associated derivations $f, g,h $ of $R$, respectively, such that $f\neq 0$ and $g\neq 0$. If  $F(x)G(y)\pm H(xy)\in Z(R)$ for all $x, y\in L$,  then $R$ commutative.
\end{cor}
\pr
Using Proposition \ref{prop-2L}, $L$ contains a nonzero ideal $I$ of $R$.  So we get  $F(x)G(y)\pm H(xy)\in Z(R)$ for all $x, y\in I$. Using Theorem \ref{th-one} and Corollary  \ref{cor-one}  we get that  $R$ is commutative.
\cqfd

Also, if we suppose $F = f$ and $G = g$ in Theorem \ref{th-one}, then we get the following corollary.

\begin{cor}
 Let $R$ be a prime ring and $I$ be a nonzero left ideal of $ R$. Let  $F,$ $G$ and $H$ be generalized derivations associated to derivations $f,$ $g$ and $h$ of $R$,  respectively,  such that $f(x)g(y)\pm H(xy)\in Z(R)$ for all $x, y\in I$. Then $R$ is commutative or $If(I)= 0$ or $Ig(I)= 0$.
\end{cor}

When $H$ in Theorem \ref{th-one} is the identity map,   we get the following result.

\begin{cor}
 Let $R$ be a  prime ring and $I$ be a nonzero left ideal of $ R$. Let  $F,$ $G$ and $H$ be generalized derivations associated to derivations $f,$ $g$ and $h$ of $R$,  respectively,  such that $F(x)G(y)\pm xy\in Z(R)$ for all $x, y\in I$. Then $R$ is commutative or $If(I)= 0$ or $Ig(I)= 0$.
\end{cor}

Now, we give our second main result. Compare it with \cite[Theorem 3]{DRR}.

\begin{thm}\label{th-tow}
 Let $R$ be a prime ring and $I$ be a nonzero left ideal of $ R$. Let  $F$  and $H$ be generalized derivations associated to derivations $f$  and $h$ of $R$,  respectively,  such that $F(x)F(y)-H(yx)\in Z(R)$ for all $x, y\in I$. If $ 0 \neq If(I) \subseteq I $  and $ 0 \neq Ih(I) \subseteq I $  , then $R$ is commutative.
 \end{thm}
\pr
 Assume that   $ I \cap Z(R)= \{0\}$. \\
 By the hypothesis, we have
\begin{equation}\label{eq-70}
    F(x)F(y)-H(yx)\in Z(R) \;\;\;\; \mbox{for  all} \;\; x,y \in I.
\end{equation}
Replacing in (\ref{eq-70}) $y$ by $yz$, where $y,z \in I$, we get,   for  all $x,y,z \in I$:
\begin{equation}\label{eq-71}
    F(x)F(y)z+F(x)yf(z)-H(yz)x-yzh(x)\in Z(R)
\end{equation}
That is
\begin{equation}\label{eq-72}
    (F(x)F(y)-H(yx))z+H(yx)z+
\end{equation}
\begin{equation*}
 F(x)yf(z)-H(yz)x-yzh(x)\in Z(R)
    \;\;\;\; \mbox{for  all} \;\; x,y,z\in I.
\end{equation*}
By the hypothesis we have
\begin{equation}\label{eq-73}
    (F(x)F(y)-H(yx))z+H(yx)z+
\end{equation}
\begin{equation*}
 F(x)yf(z)-H(yz)x-yzh(x)\in I
    \;\;\;\; \mbox{for  all} \;\; x,y,z\in I.
\end{equation*}
Since $I\cap Z(R)=0$,  we get
\begin{equation}\label{eq-74}
    (F(x)F(y)-H(yx))z+H(yx)z+
\end{equation}
\begin{equation*}
  F(x)yf(z)-H(yz)x-yzh(x)=0
    \;\;\;\; \mbox{for  all} \;\; x,y \in  I.
\end{equation*}
Replacing $z$ by $zx$, we obtain, for all $x,y,z \in  I$:
\begin{equation}\label{eq-75}
    F(x)yzf(x)-yzh(x)x-yz[x,h(x)]=0.
\end{equation}
Replacing $y$ by $py$, we get
\begin{equation}\label{eq-76}
    [F(x),p]yzf(x)=0
    \;\;\;\; \mbox{ for  all} \;\; x,y,p,z \in  I.
\end{equation}
So, using primeness of $R$, we get,  for all $x,y \in  I$:
 $[F(x),p]y=0$ or $If(x)=0$.
 This leads to
 $$
 F(x)\in Z(R) \;\;\;\mbox{or}\;\;   If(x)=0 \;\;\;\; \mbox{for  all} \;\; x \in   I.
 $$
 The fact that  a group cannot be a union of its proper subgroups with the condition     $ If(I)\neq 0 $  imply that
 $$
 F(x)\in Z(R) \;\;\; \mbox{for  all} \;\; x \in   I.
 $$
 Then, $F$ is centralizing on $ I$. Therefore,  $\dot{}R$ is commutative. Then $ I=0$,  which contradicts our hypothesis. Therefore, $I\cap Z(R)\neq \{0\}.$\\
Let  $z\in I \cap Z(R)\backslash \{0\}$ and replacing $y$ by $yz$ in our identities hypothesis, we get
\begin{equation}\label{eq80}
  F(x)yf(z)-yxh(z)\in Z(R)
   \;\;\;\mbox{for all}\;\;x, y \in I.
\end{equation}
So
\begin{equation}\label{eq81}
  [F(x)y,r]f(z)-[yx,r]h(z)=0
   \;\;\;\mbox{for all}\;\;x, y \in I;\;\; r\in R.
\end{equation}
Replacing $x$ by $xz$ in (\ref{eq81}), we get, for all $x, y \in I;\;\; r\in R$:
\begin{equation}\label{eq82}
  [F(x)(zy),r]f(z)+[xf(z)y,r]f(z)-[y(zx),r]h(z)=0.
\end{equation}
Using (\ref{eq80}) we arrive at
\begin{equation}\label{eq83}
[xf(z)y,r]f(z)=0
   \;\;\;\mbox{for all}\;\; x, y \in I ;\;\; r\in R.
\end{equation}
Replacing $x$ by $sx$ in (\ref{eq83}), where $s \in R$, we get
\begin{equation}\label{eq84}
[s,r]xf(z)yf(z)=0
   \;\;\;\mbox{for all}\;\; x, y \in I ;\;\; s \in R.
\end{equation}
The primeness of $R$ together with equation (\ref{eq84}) show that $R$
is commutative or $If(z)=0$. If $f(z)=0$,  (\ref{eq81}) becomes
\begin{equation}\label{eq85}
  [yx,r]h(z)=0
   \;\;\;\mbox{for all}\;\;x, y \in I ;\;\; r\in R.
\end{equation}
So by Lemma \ref{lem-p4} we get $h(z)=0$.\\
\\Replacing now $y$ by $tz$ in (\ref{eq-70}), where $t\in R$, we get
\begin{equation}\label{eq86}
  F(x)F(t)-H(tx)\in Z(R)
   \;\;\;\mbox{for all}\;\;x,  \in I ;\;\; t\in R.
\end{equation}
Replacing $x$ by $sz$ in (\ref{eq86}), where $s\in R$, we get
\begin{equation}
  F(s)F(t)-H(ts)\in Z(R)
   \;\;\;\mbox{for all}\;\;s,t \in R.
\end{equation}
Therefore, since $R$ is a square closed Lie ideal of $R$ itself and  by \cite[Theorem 3]{DRR} and Lemma \ref{lem-p3}, we conclude that $R$ is commutative.
\cqfd

\begin{rem}\label{rem-one}
One can see that, instead of using \cite[Theorem 3]{DRR} at the end of the proof above, we can give a direct proof following the same arguments done in \cite[Theorem 3]{DRR}.
Thus, using \cite[Theorem 1.1]{H} and Proposition \ref{prop-2L}, we conclude that   studying identities over one sided ideals can be seen until now as the ``optimal" case to study identities.
\end{rem}

We end this paper,  with some   consequences of Theorem \ref{th-tow}.

\begin{cor}
 Let $R$ be a  prime ring and $I$ be a nonzero left ideal of $ R$. Let  $F$ and $H$ be generalized derivations associated to derivations $f$ and $h$ of $R$,  respectively,  such that $F(x)F(y)+H(yx)\in Z(R)$ for all $x, y\in I$. If  $ 0 \neq If(I) \subseteq I $  and $ 0 \neq Ih(I) \subseteq I $ , then $R$ is commutative .
\end{cor}
\pr
We notice that $-H$ is a generalized derivation  of $R$   associated to the  derivation  $-h$. Hence replacing $H$ by $-H$ in Theorem \ref{th-tow}, we get $ F(x)F(y)-(-H)(yx)\in Z(R)$ for all $x, y\in I$.  That is $ F(x)F(y)+H(yx)\in Z(R)$  for all $x, y\in I$. This  implies that $R$ is commutative.
\cqfd

When we consider $F = f$  in Theorem \ref{th-tow},    we get the following result.

\begin{cor}
 Let $R$ be a  prime ring and $I$ be a nonzero left ideal of $ R$. Let  $F$  and $H$ be generalized derivations associated to derivations $f$  and $h$ of $R$,  respectively,  such that $f(x)f(y)\pm H(yx)\in Z(R)$ for all $x, y\in I$. If $ 0 \neq If(I) \subseteq I $  and $ 0 \neq Ih(I) \subseteq I $, then $R$ is commutative.
\end{cor}

Also, if we suppose $H$ to be the identity map in Theorem \ref{th-tow},  we get the following result.

\begin{cor}
 Let $R$ be a  prime ring and $I$ be a nonzero left ideal of $ R$. Let  $F$ be a generalized derivation associated to a derivation $f$  of $R$ such that $F(x)F(y)\pm yx\in Z(R)$ for all $x, y\in I$.  If $ 0 \neq If(I) \subseteq I $  and $ 0 \neq Ih(I) \subseteq I $, then $R$ is commutative.
\end{cor}

\end{document}